\documentclass[12pt,]{article}
\usepackage{pifont}
\usepackage{stmaryrd}
\usepackage[normal]{subfigure}
\usepackage{graphicx}
\usepackage{amssymb}
\usepackage{amscd}
\usepackage{amsmath}
\usepackage{amsfonts}
\usepackage{amsbsy}
\usepackage{epsfig}

\begin{document}
\date{}

\newenvironment{ar}{\begin{array}}{\end{array}}
\baselineskip20pt \addtocounter{page}{0}
\title{{\Large \bf
The maximal number of limit cycles in a family of polynomial systems\thanks{ Research supported by the Key Disciplines of Shanghai
Municipality (S30104). }}}
\author{\small{{Guanghui Xiang $^{\mbox{a}}$, Zhaoping Hu $^{\mbox{b}}$ }}
 \\{\scriptsize$^{a}$
Department of Mathematics,\ \ Shanghai Jiaotong University, Shanghai
200240, P. R. China}\\
{\scriptsize $^{b}$ Department of Mathematics,\ \ Shanghai
University, Shanghai 200444, P. R. China}} \maketitle
\begin{center}
\begin{minipage}{140mm}
\noindent{\scriptsize {{\bf Abstract.}\,\,\,The main objective of
this paper is to study the number of limit cycles in a family of
polynomial systems. Using bifurcation methods, we obtain the maximal
number of limit cycles in global bifurcation. }}

\end{minipage}
\end{center}

\section{Introduction and main results}
In the qualitative theory of real planar differential systems, the
main open problem is to determine the number and location of limit cycles. A classical
way to produce limit cycles is by perturbing a system which has a
center in such a way that limit cycles bifurcate in the perturbed
system from some of the periodic orbits of the center for the
unperturbed system. For instance, consider a planar system of the
form
\begin{eqnarray}
\begin{array}{l}
\dot{x}(t)=H_y+\varepsilon{f}(x,y,\varepsilon,a),\\
\dot{y}(t)=-H_x+\varepsilon{g}(x,y,\varepsilon,a),
\end{array}
\end{eqnarray}
where $H$, $f$, $g$ are $C^\infty$ functions in a region
$G\subset{R^2}$, $\varepsilon\in{R}$ is a small parameter and
$a\in{D}\subset{R^n}$ with $D$ compact. For $\varepsilon=0$, $(1)$
becomes a Hamiltonian system with the Hamiltonian function $H(x,y)$. Suppose
there exists a constant $H_0>0$ such that for $0<h<H_0$, the
equation $H(x,y)=h$ defines a smooth closed curve $L_h\subset{G}$
surrounding the origin and shrinking to the origin as
$h\rightarrow0$. Hence $H(0,0)=0$ and for $\varepsilon=0$ $(1)$
has a center at the origin.

Let
\begin{eqnarray}
\Phi(h,a)&=&\oint_{L_h}(gdx-fdy)_{\varepsilon=0}\\\nonumber
&=&\oint_{L_h}(H_yg+H_xf)_{\varepsilon=0}dt,
\end{eqnarray}
which is called the first order Melnikov function or Abelian
integral of $(1)$. This function plays an important role in the
study of limit cycle bifurcation. In the case that $(1)$ is a
polynomial system, a well-known problem is to study the least upper
bound of the number of zeros of $\Phi$. This is called the weakened
Hilbert $16^{th}$ problem, see [Arnold 1983; Ye, 1986].

In this paper, we first state some preliminary lemmas which can be
used to find the maximal number of limit cycles by using zeros of
$\Phi$. These lemmas are already known or easy corollaries of known
results. Then we study the global bifurcations of limit cycles for
some polynomial systems, and obtain the lower upper bound of the
number of limit cycles. This is the main part of the paper.

Now we give some lemmas. First, for Hopf bifurcation we have the
following:

\vspace{0.1cm}\textbf{Lemma 1.1([Han, 2000])} Let $H(x,y)=K(x^2+y^2)+O(|x,y|^3)$ with
$K>0$ for $(x,y)$ near the origin. Then the function $\Phi$ is of
class $C^\infty$ in $h$ at $h=0$. If
$\Phi(h,a_0)=K_1(a_0)h^{k+1}+O(h^{k+2})$, $K_1(a_0)\neq0$ for some
$a_0\in{D}$, then $(1)$ has at most $k$ limit cycles near the
origin for $|\varepsilon|+|a-a_0|$ sufficiently small.

The following lemma is well-known(see [Ye, 1986] for example).

\vspace{0.1cm}\textbf{Lemma 1.2} If
$\Phi(h,a_0)=K_2(a_0)(h-h_0)^k+O(|h-h_0|^{k+1})$,  $K_2(a_0)\neq0$
for some $a_0\in{D}$ and $h_0\in(0,H_0)$, then $(1)$ has at most
$k$ limit cycles near $L_{h_0}$ for $|\varepsilon|+|a-a_0|$
sufficiently small.

Let $L_0$ denote the origin and set
\begin{eqnarray}
S=\bigcup_{0\leq{h}<H_0}L_h.
\end{eqnarray}

It is obvious that $S$ is a simply connected open subset of the
plane. We suppose that the function $\Phi$ has the following form
\begin{eqnarray}
\Phi(h,a)=I(h)N(h,a),
\end{eqnarray}
where $I\in{C^\infty}$ for $h\in[0,H-0)$ and satisfies
\begin{eqnarray}
I(0)=0,\;\;I'(0)\neq0\;\;and\;\;I(h)\neq0\;\;for\;\;h\in(0,H_0).
\end{eqnarray}
Using above two lemmas, we can prove(see [Xiang $\&$ Han, 2004])

\vspace{0.1cm}\textbf{Lemma 1.3} Let $(4)$ and $(5)$ hold. If there exists a
positive integer $k$ such that for every $a\in{D}$ the function
$N(h,a)$ has at most $k$ zeros in $h\in[0,H_0)$ (multiplicities
taken into account), then for any given compact set $V\subset{S}$,
there exists $\varepsilon_0=\varepsilon_0(V)>0$ such that for all
$0<|\varepsilon|<\varepsilon_0$, $a\in{D}$ the system $(1)$ has at
most $k$ limit cycles in $V$.

\vspace{0.1cm}\textbf{Remark 1.1} As we known, if there exists $a_0\in{D}$ such
that the function $N(h,a)$ has exactly $k$ simple zeros
$0<h_1<\cdots<h_k<H_0$ with $N(0,a_0)\neq0$, then for any compact
set $V$ satisfying $L_{h_k}\subset{intV}$ and $V\subset{S}$, there
exists $\varepsilon_0>0$ such that for all
$0<|\varepsilon|<\varepsilon_0$, $|a-a_0|<\varepsilon_0$, $(1)$
has precisely $k$ limit cycles in $V$.

\vspace{0.1cm}\textbf{Remark 1.2} The conclusion of lemma $1.1$ and lemma $1.2$
are local with respect to both parameter $a$ and the set $S$ while
the conclusion of lemma $1.3$ is global because it holds in any
compact set of $S$ and uniformly in $a\in{D}$.

In this paper, we consider a real planar polynomial system of the
form
\begin{eqnarray}
\begin{array}{l}
\dot{x}\;=\;y(1-\alpha_1x)^{m_1}(1-\alpha_2x)^{m_2},\\
\dot{y}\;=\;-x(1-\alpha_1x)^{m_1}(1-\alpha_2x)^{m_2},
\end{array}
\end{eqnarray}
where $m_1$, $m_2$ are positive integers and $\alpha_1$, $\alpha_2$
are real constants which satisfy $\alpha_1\cdot\alpha_2\neq0$. We
shall prove that if we perturb above system by the polynomial
systems of degree $n$ we can obtain up to first order in
$\varepsilon$ at most $4([\frac{n+1}{2}]+m_1+m_2)-7$ limit cycles.

On the region
$\Omega=\{(x,y)|(1-\alpha_1x)^{m_1}(1-\alpha_2x)^{m_2}\neq0\}$, the
perturbed system by the polynomial systems of degree $n$ of
$(6)$ is equivalent to
\begin{eqnarray}
\begin{array}{l}
\dot{x}\;=\;y+\dfrac{\varepsilon}{(1-\alpha_1x)^{m_1}(1-\alpha_2x)^{m_2}}\displaystyle\sum_{0\leq{i+j}\leq{n}}a_{ij}x^iy^j,\\
\dot{y}\;=\;-x+\dfrac{\varepsilon}{(1-\alpha_1x)^{m_1}(1-\alpha_2x)^{m_2}}\displaystyle\sum_{0\leq{i+j}\leq{n}}b_{ij}x^iy^j,
\end{array}
\end{eqnarray}
where $|a_{ij}|\leq{K}$, $|b_{ij}|\leq{K}$ with $K$ a positive
constant and
$B_K=\{(a_{ij},b_{ij})|\;|a_{ij}|\leq{K},|b_{ij}|\leq{K}\}$.

Let $\Phi(h)$ denote the first order Melnikov function of $(7)$
for $0\leq{h<H_0}$,
$H_0=min(\frac{1}{\alpha_1^2},\frac{1}{\alpha_2^2})$. Then we have
the following main results.

\vspace{0.1cm}\textbf{Theorem 2.1} Suppose $\alpha_1\neq\alpha_2$. For any $K>0$
and compact set $V$ in $\Omega$, if $\Phi(h)$ is not identically
zero for $(a_{ij},b_{ij})$ varying in a compact set $D$ in $B_K$,
then there exists an $\varepsilon_0>0$ such that for
$0<|\varepsilon|<\varepsilon_0$, $(a_{ij},b_{ij})\in{D}$, the system
$(7)$ has at most $4([\frac{n+1}{2}]+m_1+m_2)-7$ limit cycles in
$V$.

\vspace{0.1cm}\textbf{Theorem 2.2} Suppose $\alpha_1=\alpha_2$. For any $K>0$ and
compact set $V$ in $\Omega$, if $\Phi(h)$ is not identically zero
for $(a_{ij},b_{ij})$ varying in a compact set $D$ in $B_K$, then
there exists an $\varepsilon_0>0$ such that for
$0<|\varepsilon|<\varepsilon_0$, $(a_{ij},b_{ij})\in{D}$, the system
$(7)$ has at most $n$ limit cycles in $V$.

\section{Proof of the theorems}

Before proving the theorems in Section 1, we give some lemmas first.

Let
\begin{eqnarray}
I_{i,j}=\oint_{L_h}\frac{x^iy^j}{(1-\alpha_1x)^{m_1}(1-\alpha_2x)^{m_2}}dt,\;\;i\geq0,j\geq0,
\end{eqnarray}
\begin{eqnarray}
I_{i,j}^{(k)}=\oint_{L_h}\frac{x^iy^j}{(1-\alpha_1x)^{k}}dt,\;\;k=1,2,3,\cdots,i\geq0,j\geq0,
\end{eqnarray}
\begin{eqnarray}
\Phi_{i,j}=a_{ij}I_{i+1,j}+b_{ij}I_{i,j+1},\;\;i\geq0,j\geq0,
\end{eqnarray}
\begin{eqnarray}
\Phi_{i,j}^{(k)}=a_{ij}I_{i+1,j}^{(k)}+b_{ij}I_{i,j+1}^{(k)},\;\;k=1,2,3,\cdots,i\geq0,j\geq0,
\end{eqnarray}
where
$$L_h:\;\;x=\sqrt{h}\sin{t},\;\;y=\sqrt{h}\cos{t}.$$
Let
$$r_1=\sqrt{1-\alpha_1^2h},\;r_2=\sqrt{1-\alpha_2^2h}.$$
The following results can be seen in the paper [Xiang $\&$ Han, 2004].

\vspace{0.1cm}\textbf{Lemma 2.1} For $m\geq1$ it holds that
\begin{eqnarray}
I_{0,0}^{(m)}=\frac{1}{r_1^{2m-1}}\sum_{j=0}^{[\frac{m-1}{2}]}C_jr_1^{2j},
\end{eqnarray}
where $C_j$$(j\geq0)$ are constants which $C_j\neq0$. $[\cdot]$
denotes the integer part function.

\vspace{0.1cm}\textbf{Lemma 2.2} For $0\leq{k}<m$, we have
\begin{eqnarray}
I_{k,0}^{(m)}&=&\sum_{j=0}^k(-1)^jC_k^jI_{0,0}^{(m-j)}\nonumber\\
&=&\frac{1}{r_1^{2m-1}}\sum_{j=0}^{[\frac{m-k-1}{2}]+k}C_jr_1^{2j},
\end{eqnarray}
and for $k\geq{m}$ we have
\begin{eqnarray}
I_{k,0}^{(m)}=\frac{1}{r_1^{2m-1}}\sum_{j=0}^{m-1}C_jr_1^{2j}+\sum_{j=0}^{[\frac{k-m}{2}]}D_jr_1^{2j},
\end{eqnarray}
where $C_j$, $D_j$ are constants.

For the function
$$\frac{x^k}{(1-\alpha_1x)^{m_1}(1-\alpha_2x)^{m_2}}$$
we have that

if $k<m_1+m_2$, there exist real constants $\tilde{A}_{k,j}$,
$\tilde{B}_{k,j}$ such that
\begin{eqnarray}
\frac{x^k}{(1-\alpha_1x)^{m_1}(1-\alpha_2x)^{m_2}}
=\sum_{j=1}^{m_1}\frac{\tilde{A}_{k,j}}{(1-\alpha_1x)^j}+\sum_{j=1}^{m_2}\frac{\tilde{B}_{k,j}}{(1-\alpha_2x)^j},
\end{eqnarray}

and if $k\geq{m_1+m_2}$, there exist real constants $A_{k,j}$,
$B_{k,j}$ and $C_{k,j}$ such that
\begin{eqnarray}
\frac{x^k}{(1-\alpha_1x)^{m_1}(1-\alpha_2x)^{m_2}}
=\sum_{j=1}^{m_1}\frac{A_{k,j}}{(1-\alpha_1x)^j}+\sum_{j=1}^{m_2}\frac{B_{k,j}}{(1-\alpha_2x)^j}+\sum_{j=0}^{k-m_1-m_2}C_{k,j}x^j.
\end{eqnarray}
Hence from the definition of $I_{k,0}$ and lemma $2.1$ for
$0\leq{k}<m_1+m_2$ we have
\begin{eqnarray}
I_{k,0}&=&\oint_{L_h}\frac{x^k}{(1-\alpha_1x)^{m_1}(1-\alpha_2x)^{m_2}}dt\nonumber\\
&=&\sum_{j=1}^{m_1}\oint_{L_h}\frac{\tilde{A}_{k,j}}{(1-\alpha_1x)^j}dt+\sum_{j=1}^{m_2}\oint_{L_h}\frac{\tilde{B}_{k,j}}{(1-\alpha_2x)^j}dt\\
&=&\frac{1}{r_1^{2m-1}}P_{m_1-1}(h)+\frac{1}{r_2^{2m-1}}P_{m_2-1}(h),\nonumber
\end{eqnarray}
and for $k\geq{m_1+m_2}$ we have
\begin{eqnarray}
I_{k,0}&=&\sum_{j=1}^{m_1}\oint_{L_h}\frac{{A}_{k,j}}{(1-\alpha_1x)^j}dt+\sum_{j=1}^{m_2}\oint_{L_h}\frac{{B}_{k,j}}{(1-\alpha_2x)^j}dt+\sum_{j=0}^{k-m_1-m_2}\oint_{L_h}C_{k,j}x^jdt\nonumber\\
&=&\frac{1}{r_1^{2m-1}}P_{m_1-1}(h)+\frac{1}{r_2^{2m-1}}P_{m_2-1}(h)+P_{[\frac{k-m_1-m_2}{2}]}(h),
\end{eqnarray}
where $P_n(h)$ denotes a polynomial of $h$ of degree $n$, and
$h=\frac{1-r_1^2}{\alpha_1^2}=\frac{1-r_2^2}{\alpha_2^2}$.

Using the definition of $L_h$ and $I_{i,j}$, we can prove
easily ([Xiang $\&$ Han, 2004])

\vspace{0.1cm}\textbf{Lemma 2.3} For $i\geq0$, $k>0$, we have
$$I_{i,2k-1}=0$$
and
$$I_{i,2k}=\sum_{j=0}^k(-1)^jC_k^jI_{i+2j,0}h^{k-j}.$$

\vspace{0.1cm}\textbf{Lemma 2.4} For $k>0$ we have
\begin{eqnarray}
\sum_{i+j=2k-1}\Phi_{ij}=\frac{1}{r_1^{2m_1-1}}P_{m_1-1+k}(h)+\frac{1}{r_2^{2m_2-1}}P_{m_2-1+k}(h)+P_{[\frac{2k-m_1-m_2}{2}]}(h)
\end{eqnarray}
and
\begin{eqnarray}
\sum_{i+j=2k}\Phi_{ij}=\frac{1}{r_1^{2m_1-1}}P_{m_1-1+k}(h)+\frac{1}{r_2^{2m_2-1}}P_{m_2-1+k}(h)+P_{[\frac{2k+1-m_1-m_2}{2}]}(h).
\end{eqnarray}

\textbf{Proof.} By the definition $\Phi_{ij}$ and lemma 2.3, we have
\begin{eqnarray*}
\sum_{i+j=2l}\Phi_{ij}&=&\sum_{i=1}^k(\Phi_{2k-2i,2i}+\Phi_{2k-2i+1,2i-1})+\Phi_{2k,0}\\
&=&\sum_{i=0}^k\tilde{a}_{2k,i}I_{2k-2i+1,2i}\\
&=&\tilde{b}_{2k,k}I_{1,0}h^k+\cdots+\tilde{b}_{2k,1}I_{2k-1,0}h+\tilde{b}_{2k,0}I_{2k+1,0}.
\end{eqnarray*}
So $(20)$ follows from $(17)$ and $(18)$. $(19)$ can be proved in the same way.

The proof is completed.

Similarly, we can prove the following formulae by using lemma $2.2$
and lemma $2.4$ of the paper [Xiang $\&$ Han, 2004].
\begin{eqnarray}
\sum_{i+j=2k-1}\Phi_{ij}^{(m)}&=&I_{0,0}^{(m)}(\tilde{b}_{2k-1,k}h^k+\cdots+\tilde{b}_{2k-1,1}h+\tilde{b}_{2k-1,0})\nonumber\\
&&+(-1)I_{0,0}^{(m-1)}(C_2^1\tilde{b}_{2k-1,k}h^k+\cdots+C_{2k-2}^1\tilde{b}_{2k-1,1}h+C_{2k}^1\tilde{b}_{2k-1,0})\nonumber\\
&&\vdots\nonumber\\
&&+(-1)^{m-1}I_{0,0}^{(1)}(C_{2[\frac{m}{2}]}^{m-1}\tilde{b}_{2k,k-[\frac{m}{2}]}h^{k-[\frac{m}{2}]}+\cdots+C_{2k-2}^{m-1}\tilde{b}_{2k-1,1}h+C_{2k}^{m-1}\tilde{b}_{2k-1,0})\nonumber\\
&&\nonumber\\
&&+(-1)^mC_{2[\frac{m+1}{2}]-1}^{m-1}\Big(K_0\tilde{b}_{2k-1,k-[\frac{m+1}{2}]}+\cdots+K_{2k-2[\frac{m+1}{2}]}\tilde{b}_{2k-1,0}\Big)h^{k-[\frac{m+1}{2}]}\nonumber\\
&&\vdots\nonumber\\
&&+(-1)^mC_{2k-3}^{m-1}\Big(K_0\tilde{b}_{2k-1,1}+K_2\tilde{b}_{2k-1,0}\Big)h+(-1)^mC_{2k-1}^{m-1}K_0\tilde{b}_{2k-1,0}
\end{eqnarray}
and
\begin{eqnarray}
\sum_{i+j=2k}\Phi_{ij}^{(m)}&=&I_{0,0}^{(m)}(\tilde{b}_{2k,k}h^k+\cdots+\tilde{b}_{2k,1}h+\tilde{b}_{2k,0})\nonumber\\
&&+(-1)I_{0,0}^{(m-1)}(C_2^1\tilde{b}_{2k,k}h^k+\cdots+C_{2k-1}^1\tilde{b}_{2k,1}h+C_{2k+1}^1\tilde{b}_{2k,0})\nonumber\\
&&\vdots\nonumber\\
&&+(-1)^{m-1}I_{0,0}^{(1)}(C_{2[\frac{m+1}{2}]-1}^{m-1}\tilde{b}_{2k,k-[\frac{m-1}{2}]}h^{k-[\frac{m}{2}]}+\cdots+C_{2k-1}^{m-1}\tilde{b}_{2k,1}h+C_{2k+1}^{m-1}\tilde{b}_{2k,0})\nonumber\\
&&\nonumber\\
&&+(-1)^mC_{2[\frac{m}{2}]}^{m-1}\Big(K_0\tilde{b}_{2k,k-[\frac{m}{2}]}+\cdots+K_{2k-2[\frac{m}{2}]}\tilde{b}_{2k,0}\Big)h^{k-[\frac{m}{2}]}\nonumber\\
&&\vdots\nonumber\\
&&+(-1)^mC_{2k-2}^{m-1}\Big(K_0\tilde{b}_{2k,1}+K_2\tilde{b}_{2k,0}\Big)h+(-1)^mC_{2k}^{m-1}K_0\tilde{b}_{2k,0}
\end{eqnarray}

Now we are in position to prove the main results.

\vspace{0.1cm}\textbf{Proof of Theorem 2.1} In the following we suppose $n=2s$
first. In this case, by $(2)$ the Melnikov function $\Phi(h)$ of
system $(7)$ has the following form
\begin{eqnarray}
\Phi(h)&=&\oint_{L_h}\frac{1}{(1-\alpha_1x)^{m_1}(1-\alpha_2x)^{m_2}}\sum_{0\leq{i+j}\leq2s}(a_{ij}x^{i+1}y^j+b_{ij}x^iy^{j+1})dt\nonumber\\
&=&\sum_{0\leq{i+j}\leq2s}\Phi_{ij}\\
&=&\sum_{k=1}^s(\sum_{i+j=2k-1}\Phi_{ij}+\sum_{i+j=2k}\Phi_{ij})+\Phi_{00}.\nonumber
\end{eqnarray}
From $(19)$ and $(20)$, $(23)$ becomes
\begin{eqnarray}
\Phi(h)&=&\sum_{k=0}^s\Big(\frac{1}{r_1^{2m_1-1}}P_{m_1-1+k}(h)+\frac{1}{r_2^{2m_2-1}}P_{m_2-1+k}(h)+P_{[\frac{2k+1-m_1-m_2}{2}]}(h)\Big)\nonumber\\
&=&\frac{1}{r_1^{2m_1-1}}P_{m_1-1+s}(h)+\frac{1}{r_2^{2m_2-1}}P_{m_2-1+s}(h)+P_{[\frac{2s+1-m_1-m_2}{2}]}(h),
\end{eqnarray}
where $P_k(h)$ is a polynomial of $h$ of degree $k$.

Obviously all the zeros of $(24)$ satisfy
$$\Big[\frac{1}{r_1^{2m_1-1}}P_{m_1-1+s}(h)+\frac{1}{r_2^{2m_2-1}}P_{m_2-1+s}(h)\Big]^2=\Big[P_{[\frac{2s+1-m_1-m_2}{2}]}(h)\Big]^2.$$
Further the above formula becomes
$$\sqrt{(1-\alpha_1^2h)(1-\alpha_2^2h)}P_{2s+2(m_1+m_2)-4}(h)=Q_{2s+2(m_1+m_2-3)}(h),$$
where $P_{2s+2(m_1+m_2)-4}(h)$ and $Q_{2s+2(m_1+m_2-3)}(h)$ are two
real coefficient polynomials of $h$ of degree $2s+2(m_1+m_2)-4$ and
$2s+2(m_1+m_2)-3$ respectively. Hence the number of zeros of
$\Phi(h)$ are not large than $4s+4(m_1+m_2)-6$.

For the case of $n=2s-1$, similarly we can prove that the number of
zeros of $\Phi(h)$ are not large than $4s+4(m_1+m_2)-6$.

Notice that $\Phi(h)=0$ at $h=0$ in $(23)$. From lemma $1.3$, we
know that there exists an $\varepsilon_0>0$ such that when
$0<|\varepsilon|<\varepsilon_0$, $a=(a_{ij},b_{ij})$ which satisfy
$|a_{ij}|\leq{K}$, $|b_{ij}|\leq{K}$, the system $(9)$ has at most
$4\Big(m_1+m_2+\frac{n+1}{2}\Big)-7$ limit cycles. The proof is
completed.

\vspace{0.1cm}\textbf{Proof of theorem 2.2} We suppose $\alpha_1=\alpha_2$ and
$m_1+m_2=m$ in $(7)$.

For the case of $n=2s$, from $(21)$ and $(22)$ the
Melnikov function $\Phi(h)$ of system $(7)$ has the following
form
\begin{eqnarray*}
\Phi(h)&=&I_{(0,0)}^{(m)}(b_s^{(m)}h^s+\cdots+b_1^{(m)}h+b_0^{(m)})\\
&&+\cdots\\
&&+I_{0,0}^{(1)}(b_{s-[\frac{m-1}{2}]}^{(1)}h^{s-[\frac{m-1}{2}]}+\cdots+b_1^{(1)}h+b_0^{(1)})\\
&&+(B_{s-[\frac{m}{2}]}h^{s-[\frac{m}{2}]}+\cdots+B_1h+B_0),
\end{eqnarray*}
where $b_j^{(i)}$, $B_j$$(1\leq{i}\leq{m},j\geq0)$ are linear
combinations of $a_{ij}$, $b_{ij}$ with $0\leq{i+j}\leq2s$. Let
$\sqrt{1-\alpha_1^2h}=r$, $0<r<1$. And from $(12)$ and $(13)$, the
above formula becomes
\begin{eqnarray*}
\Phi(h)&=&\frac{1}{r^{2m-1}}(c_{2s+m}r^{2s+m}+c_{2s+m-1}r^{2s+m-1}+\cdots+c_{2m-1}r^{2m-1}\\
&&+c_{2m-2}r^{2m-2}+c_{2m-4}r^{2m-4}+\cdots+c_2r^2+c_0)\\
&=&\frac{1}{r^{2m-1}}P_{2s+m}(r),
\end{eqnarray*}

 where $P_{2s+m}(r)$ is a polynomial of $r$ of degree $2s+m$
and $P_{2s+m}(r)=0$ at $r=1$. Notice that the polynomial
$P_{2s+m}(r)$ has only $2s+2$ items. By Rolle theorem $P_{2s+m}(r)$
has at most $2s+1$ positive zeros. So the polynomial
$\frac{P_{2s+m}(r)}{1-r}$ has at most $2s$ positive zeros. From
lemma $1.3$, we know that there exists an $\varepsilon_0>0$ such
that when $0<|\varepsilon|<\varepsilon_0$, $a=(a_{ij},b_{ij})$ which
satisfy $|a_{ij}|\leq{K}$, $|b_{ij}|\leq{K}$ the system $(7)$ has
at most $2s$ limit cycles.

For the case of $n=2s-1$ we can prove the theorem in the same way.
The proof is completed.

\vspace{0.1cm}\textbf{Remark 2.1} In fact, for the system
\begin{eqnarray*}
\begin{array}{l}
\dot{x}\;=\;y(1-\alpha_1x)^{m_1}(1-\alpha_2x)^{m_2}\cdots(1-\alpha_kx)^{m_k},\\
\dot{y}\;=\;-x(1-\alpha_1x)^{m_1}(1-\alpha_2x)^{m_2}\cdots(1-\alpha_kx)^{m_k},
\end{array}
\end{eqnarray*}
where $m_1,m_2,\cdots,m_k$ are positive integers and
$\alpha_1,\alpha_2,\cdots,\alpha_k$ are real constants which satisfy
$\alpha_1\cdot\alpha_2\cdot\cdots\cdot\alpha_k\neq0$. Using the same
way we can prove that if we perturb the above system inside the
polynomial systems of degree $n$ we can obtain up to first order in
$\varepsilon$ at most
$2^{k}\Big([n\frac{n+1}{2}]+\sum_{j=1}^{k}m_j-k\Big)+2^{k-1}(k-1)-1$
limit cycles.

\vspace*{1cm}

\noindent {\bf\Large References}

\begin{description}

\item[] Arnold V.[1983] ``Geometrical Methods in the Theory of Ordinary Differential
Equations,`` {\em Springer-Verlag, New York}.

\item[] Cima A., Gasull A. $\&$ Manosas F.[1995] ``Cyclicity of a family of
vectorfields," {\em J. Math. Anal. Appl.} {\bf 196}, 921-937.

\item[] Han M. [1999] ``On cyclicity of planar systems in Hopf and Poincare
bifurcations," {\em Dynamical Systems Proceedings of the International
Conference in Honor of Professor Liao Shantao}, 66-74.

\item[] Han M. [2000] ``On hopf cyclicity of planar systems," {\em J. Math. Anal.
Appl} {\bf 245}, 404-422.

\item[] Han M. [2001] ``The Hopf cyclicity of Lineared systems," {\em Appl. Math.
Letters} {\bf 14}, 183-188.

\item[] Li C., Libre J. $\&$ Zhang Z. [1995] ``Weak focus, limit cycles and
bifurcations for bounded quadratic systems," {\em J. Diff. Eqs.} {\bf 115}, 193-223.
\item[] Li C., Li W., Libre J. $\&$ Zhang Z. [2000] ``Linear estimate for the number of
zeros of Abelian integrals for quadratic isochronous centres," {\em Nonlinearity} {\bf 13}, 1775-1800.

\item[] Novikov D. $\&$ Yakovenko S.[1995] ``Simple exponential estimate for the
number of real zeros of complete Abelian integrals," {\em Ann. Inst. Fourier.} {\bf 45}, 897-927.

\item[] Varchenko A. [1984] ``Estimate of the number of zeros of Abelian
integrals depending on parameters and limit cycles," {\em Funt. Anal. Appl.} {\bf 18}, 98-108.

\item[] Xiang G. $\&$ Han M. [2004] ``Global bifurcation of limit cycles in a fammily
of multiparameter systems," {\em Int. J. Bifurcation and Chaos} {\bf 14}, 3325-3335.

\item[] Ye Y. [1986] ``Theory of limit cycles," {\em Transl Math Monographs, Amer Math Soc} {\bf 66}.

\end{description}

\end{document}